\documentclass[12pt]{amsart}

\usepackage{latexsym}
\usepackage[centertags]{amsmath}
\usepackage{amsfonts}
\usepackage{amssymb}
\usepackage{amsthm}
\usepackage{newlfont}
\usepackage[square,sort]{natbib}

\def\mytopsep{3mm}

\newtheoremstyle{myplain}{\mytopsep}{\mytopsep}{\itshape}{0pt}{\bfseries}{.}{3mm}{}
\newtheoremstyle{mydefinition}{\mytopsep}{\mytopsep}{\normalfont}{0pt}{\bfseries}{.}{3mm}{}

\newtheoremstyle{myremark}{\mytopsep}{\mytopsep}{\normalfont}{0pt}{\bfseries}{.}{3mm}{}

\theoremstyle{myplain}

\newtheorem{thm}{Theorem}[section]

\newtheorem{lem}[thm]{Lemma}
\newtheorem{prop}[thm]{Proposition}

\theoremstyle{mydefinition}
\newtheorem{dfn}[thm]{Definition}
\theoremstyle{myremark}
\newtheorem{rem}[thm]{Remark}
\newtheorem{exa}[thm]{Example}

\allowdisplaybreaks[1]

\makeatletter\@addtoreset{equation}{section}\makeatother

\usepackage{maple2e}

\DefineParaStyle{maplegroup} \DefineParaStyle{Author}
\DefineParaStyle{subtitle} \DefineParaStyle{Bullet Item}
\DefineParaStyle{Dash Item} \DefineParaStyle{Diagnostic}
\DefineParaStyle{Error} \DefineParaStyle{Heading 1}
\DefineParaStyle{Heading 2} \DefineParaStyle{Heading 3}
\DefineParaStyle{Heading 4} \DefineParaStyle{Normal}
\DefineParaStyle{Text Output} \DefineParaStyle{Title}
\DefineParaStyle{Warning} \DefineCharStyle{2D Math}
\DefineCharStyle{Help Heading} \DefineCharStyle{Maple Input}
\DefineCharStyle{FixedWidth} \DefineCharStyle{endFixedWidth}

\def\Oge{\mathop{\Omega}_\ge}
\def\Oeq{\mathop{\Omega}_=}

\def\supp{\mathop{\mbox{supp}}}
\def\ord{\mathrm{ord}}
\def\mb{\mathbf}

\def\NN{\mathbb{N}}
\def\tt{t^{-1}}

\def\PP{\mathbb{P}}

\def\CC{\mathbb{C}}
\def\ZZ{\mathbb{Z}}

\def\res{\mathop{\mathrm{Res}}}
\def\ct{\mathop{\mathrm{CT}}}
\def\PT{\mathop{\mathrm{PT}}}

\def\pt{\mathop{\mathrm{PT}}}

\def\poly{\mbox{\sf Poly}}
\def\frr{\mbox{\sf Frac}}

\renewcommand{\ll}{\langle\!\langle}
\renewcommand{\gg}{\rangle\!\rangle}

\begin{document}

\title{A Fast Algorithm for MacMahon's Partition Analysis}

\author{Guoce Xin}
\address{Department
of Mathematics\\
Brandeis University\\
Waltham MA 02454-9110} \email{maxima@brandeis.edu}

\date{August 9, 2004}
\begin{abstract}

This paper deals with evaluating constant terms of a special class
of rational functions, the Elliott-rational functions. The
constant term of such a function can be read off immediately from
its partial fraction decomposition. We combine the theory of
iterated Laurent series and a new algorithm for partial fraction
decompositions to obtain a fast algorithm for MacMahon's Omega
calculus, which (partially) avoids the ``run-time explosion"
problem when eliminating several variables. We discuss the
efficiency of our algorithm by investigating problems studied by
Andrews and his coauthors; our running time is much less than that
of their Omega package.

\end{abstract} \maketitle

\section{Introduction}

Zeilberger \cite{zeil} proved a conjecture of Chan et al.
\cite{chan} by proving an identity {equivalent to}
\begin{equation}\label{e-zeil}
\ct_{x_1}\cdots \ct_{x_n} \frac{1}{\prod_{i=1}^n (1-x_i)}
\frac{1}{\prod_{i<j} (x_i-x_j)} =C_1\cdots C_{n-1},
\end{equation}
where $C_k$'s are the Catalan numbers.

This identity should be interpreted as taking \emph{{iterated
constant terms}} \cite{welleda}; i.e., in applying $\ct_{x_n}$ to
the displayed rational function, we expand it as a Laurent series
in $x_n$; the result is still a rational function and we can apply
$\ct_{x_{n-1}}$, \dots, $\ct_{x_1}$ to it iteratively.

The idea behind the above treatment is to give a proper series
expansion of $1/(x_i-x_j)$ for every $i$ and $j$, so that all of
the expansions are compatible. Once we have determined the
relations between the $x$'s, there is no confusion about their
series expansion. For instance, we can let $1>x_1>\cdots >x_n$.
For the particular rational function in equation \eqref{e-zeil},
which is symmetric in the $x$'s, there is no confusion after a
total ordering on the $x$'s is given.

Here we present a slightly different, but more efficient,
approach, by means of applying the theory of the field of
\emph{iterated Laurent series}. We first treat the rational
function in question as an {iterated Laurent series}, by which we
mean we expand it as a Laurent series in $x_n$, then a Laurent
series in $x_{n-1}$, and so on. Then we take the constant term.
This idea led to the study of the field of iterated Laurent series
in \cite[Ch. 2]{xinthesis}, which applies to MacMahon's Partition
Analysis.

MacMahon's Partition Analysis is suited for solving problems of
counting solutions to linear Diophantine equations and
inequalities. Using MacMahon's approach, problems such as counting
lattice points in a convex polytope, counting integral solutions
to a system of linear Diophantine equations, and computing Ehrhart
quasi-polynomials, become evaluations of the constant term of an
\emph{Elliott-rational function}: a rational function whose
denominator has only factors of the form $A-B$, where $A$ and $B$
are both monomials. An example of such is the rational function in
\eqref{e-zeil}.

MacMahon's technique has been restudied by Andrews et. al. using
computer algebra in a series papers [1--9].
New algorithms have been found and computer programs such as the
Omega package have been developed.

The constant term (in one variable) of an Elliott-rational
function can be read off immediately if its partial fraction
decomposition is given. However, the coefficients of a rational
function must lie in a field to guarantee the existence of its
partial fraction decompositions, and the classical algorithm for
partial fraction decomposition is rather slow because the
coefficients contain many other variables. The above two problems
are solved by applying the theory of iterated Laurent series and a
new algorithm for partial fraction decompositions in
\cite{xinparfrac}.

In section 2, we give the basic theory of iterated Laurent series.
The fundamental structure theorem tells us when a formal Laurent
series is an iterated Laurent series. In section 3, we introduce
MacMahon's partition analysis. In section 4, we develop an
efficient algorithm for MacMahon's partition analysis by combining
the theory of iterated Laurent series and a new algorithm for
partial fraction decompositions. The theory of iterated Laurent
series is crucial in avoiding the ``run-time explosion" problem
\cite[p. 9]{andrews7} when eliminating several variables. In
section 5, we use our Maple package to test the efficiency of our
algorithm. We investigate problems related to $k$-gons,
generalized Putnam problems, and magic squares
\cite{andrews5,george6,andrews7,andrews9}. The known formulas are
obtained within seconds, and several new formulas are produced in
minutes. Finally in section 6, we point out several ways to
accelerate the computer program. There are also ways to make the
computation easier that are hard to implement on the computer. As
an example, we give a simple proof of the formula for $k$-gon
partitions in \cite{andrews9}.

\section{The Field of Iterated Laurent Series}

By a formal Laurent series in $x_1,\dots ,x_n$, we mean a series
that can be written in the form
$$\sum_{i_1=-\infty}^{\infty}\cdots \sum_{i_n=-\infty}^{\infty} a_{i_1\ldots i_n}x_1^{i_1}
\cdots x_n^{i_n},$$ where $a_{i_1\ldots i_n}$ are elements in a
field $K$. For formal Laurent series, the definition of the
constant term operator is clear:

\begin{dfn}[Natural Definition]\label{dfn-natural}
The operator $\ct_{x_j}$ acts on a formal series in $x_1,\dots
,x_n$ with coefficients $a_{i_1,\dots ,i_n}$ in $K$ by
$$\ct_{x_j} \sum_{(i_1,\dots ,i_n)\in \ZZ^n}
a_{i_1,\dots ,i_n} x_1^{i_1} \cdots x_n^{i_n} = \sum_{(i_1,\dots
,i_n)\in \ZZ^n, i_j=0} a_{i_1,\dots ,i_n} x_1^{i_1} \cdots
x_n^{i_n}.$$
\end{dfn}

The simplest way to apply the natural definition would be to work
with all formal series $\sum_{i_1,\dots ,i_n} a_{i_1,\dots ,i_n}
x_1^{i_1} \cdots x_n^{i_n}$, where $(i_1,\dots ,i_n)$ ranges over
all elements of $\ZZ^n$. Unfortunately, they do not form a ring.
Therefore we usually work in a ring, such as the ring of Laurent
series $K((x_1,\dots ,x_n))$: formal series of monomials where the
exponents of the variables are bounded from below. But we need a
larger ring or even a field that includes all rational functions,
because many constant term evaluation problems involves rational
functions.

Let $K$ be a field. We define $K\ll x_1\gg$ to be the field of
Laurent series $K((x_1))$, and define the {\em field of iterated
Laurent series} $K\ll x_1,\dots ,x_n\gg$ inductively to be $K\ll
x_1,\dots ,x_{n-1}\gg((x_n))$, which is the field of Laurent
series in $x_n$ with coefficients in $K\ll x_1,\dots ,x_{n-1}\gg$.
Thus an iterated Laurent series is first regarded as a Laurent
series in $x_n$, then a Laurent series in $x_{n-1}$, and so on. An
iterated Laurent series obviously has a unique formal Laurent
series expansion. However, it is not obvious which formal series
are in $K\ll x_1,\dots ,x_n\gg$. The fundamental structure theorem
solves this problem nicely.

We define a total ordering $\preceq$ on monomials by representing
$x_1^{i_1}\cdots x_n^{i_n}$ by $(i_1,\dots, i_n)\in \ZZ^n$, where
$\ZZ^n$ is ordered reverse lexicographically. So $x_i^s\prec x_j$
for all $i<j$ and $s\in \ZZ$. We define the \emph{support} of a
formal Laurent series by
$$\supp\ \sum_{(i_1,\dots ,i_n)\in \ZZ^n}
a_{i_1,\dots ,i_n} x_1^{i_1} \cdots x_n^{i_n} := \{\, (i_1,\dots
i_n)\mid a_{i_1,\dots ,i_n} \ne 0 \,\}.
$$
Recall that a totally ordered set $S$ is {\em well-ordered} if
each nonempty subset of $S$ contains a minimal element.

\begin{thm}[Fundamental Structure] \label{p-wellorder} {A formal series in
$x_1,\dots ,x_n$ belongs
to \\
$K\ll x_1,\dots ,x_{n}\gg$ if and only if it has a well-ordered
support.}
\end{thm}
The proof of this theorem is omitted. For details, see
\cite[Proposition 2.1.2]{xinthesis}. The result gives us an
overview about when a formal Laurent series is an iterated Laurent
series.

The fundamental structure theorem, together with the simple and
useful fact that \emph{any subset of a well-ordered set is
well-ordered}, justify the application of the natural definition
in $K\ll x_1,\dots ,x_n \gg$ because of the following three
properties:
\begin{enumerate}
\item[$P1$.] $\ct_{x_i}:K\ll x_1,\dots ,x_n \gg
\rightarrow K\ll x_1,\dots, \hat{x}_i , \dots,x_n \gg.$ This
property is necessary to make the natural definition applicable.

\item[$P2$.] $\ct_{x_k}\sum_i F_i=\sum_i \ct_{x_k} F_i$.
This property is the key to converting many problems into simple
algebraic computations.

\item[$P3$.] $\ct_{x_i}\ct_{x_j}F=\ct_{x_j}\ct_{x_i}F$.
This property may significantly simplify the constant term
evaluations.
\end{enumerate}

We define the \emph{order} $\ord(f)$ of an iterated Laurent series
$f$ to be the minimum of its support, which is well-ordered by the
fundamental structure theorem. We have the following composition
law.

\begin{prop}[Composition Law]
Suppose that $f$ belongs to $K\ll x_1,\dots ,x_n\gg$ and
$\ord(f)>\ord(1)$. Then for any $b_i\in K$ for all $i$,
$$\sum_{i=0}^\infty b_i f^i$$
is well defined and belong to $K\ll x_1,\dots ,x_n\gg$, in the
sense that all of its coefficients are finite sum of nonzero
elements in $K$.
\end{prop}

This result is a consequence of a general result for
Malcev-Neumann series \cite[Theorem 3.1.7]{xinthesis}. As a
consequence, the series expansion of $1/(1-f)$ for
$\ord(f)>\ord(1)$ is just $1+f+f^2+\cdots$. More generally, for
two iterated Laurent series $A$ and $ B$ with $\ord(A)<\ord(B)$,
the expansion of $1/(A-B)$ is
$$ \frac{1}{A-B}=\frac{1}{A}\frac{1}{1-B/A}= \sum_{k\ge 0} B^k/A^{k+1}.$$
For instance, in $K\ll x_1,x_2, x_3 \gg$, we have
$$ \frac{1}{x_1^2x_2^4-x_3}= \sum_{k\ge 0} x_3^k/(x_1^2x_2^4)^{k+1}.$$

In the field $K\ll x_1,\dots x_n\gg$, we define a total ordering
on the variables, which produces a total ordering on its group of
monomials. This total ordering plays a central role in series
expansion. By thinking of iterated Laurent series as numbers,
$\ord(f)>\ord(1)$ means that $f$ is much smaller than $1$, or
$f=o(1)$. Similarly $\ord(B)>\ord(A)$ means that $B$ is much
smaller than $A$, or $B=o(A)$.

The analogous situation for complex variables would be informally
written as $1>\!\!> x_1>\!\!> \cdots
>\!\!> x_n$ when expanding rational functions into Laurent series,
where $>\!\!> $ means ``much greater". See \cite{wilson} and
\cite[p. 231]{stanley-rec}.

The the following three computational rules are frequently used in
constant term evaluations. Let $F,G\in K\ll x_1,\dots ,x_n\gg$.

1. Linearity: $\ct_{x_i} (aF+bG)= a \ct_{x_i}F+b\ct_{x_i} G$, if
$a$ and $b$ are independent of $x_i$.

2. If $F$ can be written as $\sum_{k\ge 0} a_k x_i^k$, then
$\displaystyle \ct_{x_i} F = \left. F \right|_{x_i=0}.$

3. $\displaystyle\res_{x_i} \frac{\partial}{\partial x_i} F =0$.

\begin{rem}
Depending on the working field, rational functions
$Q(x_1,x_2,\ldots,x_m)$ may have as many as $m!$ different
expansions. More precisely, if $\sigma$ is a permutation of $[m]$,
then $Q(\mb{x})$ will have a unique expansion in $K\ll
x_{\sigma_1},x_{\sigma_2},\ldots ,x_{\sigma_m}\gg $. The
expansions of $Q(\mb{x})$ for different $\sigma$ are usually
different. So we need to specify the working field whenever a
reciprocal comes into account.
\end{rem}

Iterated Laurent series is to obtained by defining a total
ordering on its variables (this idea is not new, e.g.,
\cite{stanley-rec, wilson}). In fact, it is a special kind of
Malcev-Neumann series, which has been studied in \cite{xinthesis},
and has applications to MacMahon's partition analysis.

\section{MacMahon's Partition Analysis}

MacMahon's Partition Analysis is used for counting the solutions
to a system of linear Diophantine equations and inequalities, and
 the number of lattice points in a convex polytope.
Such problems can be converted into evaluating the constant terms
of certain \emph{Elliott-rational functions}. This conversion has
been known as MacMahon's partition analysis, and has been given a
new life by Andrews et al. in a series of papers [1--9].

\begin{dfn}
An {\em Elliott-rational} function is a rational function that can
be written in such a way that its denominator can be factored into
the products of one monomial minus another, with the $0$ monomial
allowed.
\end{dfn}

In the one-variable case, this concept reduces to the generating
function of a quasi-polynomial.

MacMahon's idea was to introduce new variables
$\lambda_1,\lambda_2,\dots $ to replace linear constraints. For
example, suppose we want to count the nonnegative integral
solutions to the linear equation $2a_1-3a_2+a_3+2=0$. We can
compute the generating function of such solutions as the
following:

$$\sum_{a_1,a_2,a_3\ge 0 \atop 2a_1-3a_2+a_3+2=0} x_1^{a_1}x_2^{a_2}x_3^{a_3}
=\sum_{a_1,a_2,a_3 \ge 0} \ct_\lambda \lambda^{2a_1-3a_2+a_3+2}
x_1^{a_1}x_2^{a_2}x_3^{a_3}.$$ Now apply the formula for the sum
of a geometric series. It becomes
$$\ct_\lambda \frac{\lambda^2}{(1-\lambda^2 x_1)(1-\lambda^{-3}x_2)(1-\lambda x_3)}.$$
The above expression is a power series in $x_i$ but not in
$\lambda$.

It is clear that if there are $r$ linear equations, we can compute
their solutions by introducing $r$ variables $\Lambda$, short for
$\lambda_1,\dots ,\lambda_r$. Thus counting solutions of a system
of linear Diophantine equations can be converted into evaluating
the constant term of an Elliott-rational function.

\begin{thm}\label{t-elliott}
If $F$ is Elliott-rational, then the constant terms of $F$ are
still Elliott-rational.
\end{thm}

This result follows from  ``The method of Elliott" (see \cite[p.
111--114]{mac}) developed from the following identity. Note that
we have not specified the working field yet.

\begin{lem}[Elliott Reduction Identity]
For positive integers $j$ and $k$,
$$\frac{1}{(1-x\lambda^j)(1-y\lambda^{-k})}=\frac{1}{1-xy\lambda^{j-k}}
\left(\frac{1}{1-x\lambda^j}+\frac{1}{1-y\lambda^{-k}}-1\right).$$
\end{lem}

Elliott's argument is that after finitely many applications of the
above identity to an Elliott-rational function, we will get a sum
of rational functions, in which every denominator has either all
factors of the form $1-x\lambda^i$, or all factors of the form
$1-y/\lambda^i$. Now taking the constant term of each summand is
easy.

Theorem \ref{t-elliott} reduces the evaluation of $\ct_{\Lambda}
F$ to the univariate case $\ct_{\lambda}F$ by iteration.
Unfortunately, the Elliott reduction algorithm is not efficient in
practice. Other algorithms have been developed, and computer
programs have been set up, such as the ``Omega" package
\cite{george6}. But we can do much better by the partial fraction
method and working in a field of iterated Laurent series.

Before going further, let us review some of the work in
\cite{george6}.
 The key ingredient in
their argument is MacMahon's Omega operator $\Omega_\ge$, which is
defined by:
$$\Oge \sum_{s_1=-\infty}^\infty \cdots \sum_{s_r=-\infty}^{\infty} A_{s_1,\dots ,s_r}
\lambda_1^{s_1} \cdots \lambda_r^{s_r} := \sum_{s_1=0}^\infty
\cdots \sum_{s_r=0}^{\infty} A_{s_1,\dots ,s_r},$$ where the
domain of the $A_{s_1,\dots ,s_r}$ is the field of rational
functions over $\CC$ in several complex variables and $\lambda_i$
are restricted to a neighborhood of the circle $|\lambda_i|=1.$ In
addition, the $A_{s_1,\dots ,s_r}$ are required to be such that
any of the $2^{r}-1$ sums
$$\sum_{s_{i_1}=0}^\infty \cdots \sum_{s_{i_j}=0}^{\infty} A_{s_{i_1},\dots
,s_{i_j}}$$
 is absolute convergent within the domain of the definition of $A_{s_1,\dots ,s_r}$.

Another operator $\Oeq$ is given by
$$\Oeq \sum_{s_1=-\infty}^\infty \cdots \sum_{s_r=-\infty}^{\infty} A_{s_1,\dots ,s_r}
\lambda_1^{s_1} \cdots \lambda_r^{s_r} := A_{0,\dots ,0}.$$

Andrews et al. emphasized in \cite{george6} that it is essential
to treat everything analytically rather than formally because the
method relies on unique Laurent series representations of rational
functions.

It is not hard to see that their definition always works if we are
working in a ring such as the ring of formal power series in
$\mb{x}$ with coefficients Laurent polynomials in $\Lambda$, where
$\mb{x}$ is short for $x_1,\dots, x_n$ and $\Lambda$ is short for
$\lambda_1,\dots ,\lambda_r$. In fact, this approach was used by
Han in \cite{han}.

By Theorem \ref{t-elliott}, it suffices to consider the case of
$r=1$, since the general case can be done by iteration. In the
previous work by Andrews et al. and by Han, the problem was
reduced to evaluating the constant term (with respect to
$\lambda$) of a rational function of the form
\begin{align}\label{e-2-mac-g}
\frac{\lambda^k}{\prod_{1\le i\le m}
(1-\lambda^{j_i}x_i)\prod_{1\le i\le n} (1-y_i/\lambda^{k_i})}.
\end{align}
This treatment has assumed the obvious geometric expansion:
$$\frac{1}{1-\lambda^{j_i}x_i}=\sum_{s=0}^\infty \lambda^{sj_i}x_i^s\quad
\text{and}\quad \frac{1}{1-y_i/\lambda^{k_i}}=\sum_{s=0}^\infty
\lambda^{-sk_i}y_i^s.$$ In other words, for each factor $f$ in the
denominator, $f$ has positive powers in $\lambda$ indicates that
the series expansion of $1/f$ contains only nonnegative powers in
$\lambda$; and $f$ has negative powers in $\lambda$ indicates that
$1/f$ contains only nonpositive powers in $\lambda$. In our
approach, these indications are dropped off after defining a total
ordering.

We find it better to do this kind of work in a certain field of
iterated Laurent series, because in such a field, we can use the
theory of partial fraction decompositions in $K(\lambda)$ for any
field $K$ and any variable $\lambda$. We illustrate this idea by
solving a problem in \cite[p. 252]{george6} with the partial
fraction method.

\vspace{3mm} \noindent {\bf Problem.} Find all nonnegative integer
solutions $a,b$ to the inequality $2a\ge 3b$.

\begin{proof}[Solution] First of all, using geometric series summations  we translate
the problem into a form which MacMahon calls the {\em crude
generating function}, namely
$$f(x,y):=\sum_{a,b\ge 0, 2a-3b\ge 0} x^a y^b =\Oge \sum_{a,b\ge 0} \lambda^{2a-3b} x^ay^b
=\Oge \frac{1}{(1-\lambda^2 x)(1-\lambda^{-3}y)},$$ where
everything is regarded as a power series in $x$ and $y$ but not in
$\lambda$.

Now by converting into partial fractions in $\lambda$, we have
$$\frac{1}{(1-\lambda^2 x)(1-\lambda^{-3}y)}=
\frac{y(1+\lambda x^2y+ \lambda^2 x)}{(1-x^3y^2)(\lambda^3-y)}+
\frac{1+\lambda x^2y}{ (1-x^3y^2)(1-\lambda^2 x)}.$$ Where the
right-hand side of the above equation is expanded as a power
series in $x$ and $y$, the second term contains only nonnegative
powers in $\lambda$, and the first term,
$$\frac{y(1+\lambda x^2y+ \lambda^2 x)}{(1-x^3y^2)(a^3-y)}=\frac{y}{1-x^3y^2}
\frac{\lambda^{-3}+\lambda^{-2}x^2y +\lambda^{-1} x}{
1-\lambda^{-3} y}$$ contains only negative powers in $\lambda$.
Thus by setting $\lambda=1$ in the second term, we obtain
$$f(x,y)=\frac{1+x^2y}{
(1-x^3y^2)(1- x)}.$$ By a geometric series expansion, it is easy
to deduce that
$$\{\, (a,b)\in \NN^2: 2a\ge 3b \,\} = \{\, (m+n+\lceil n/2 \rceil ,n): (m,n)\in \NN^2\,\}.$$
\end{proof}

In solving the above problem, we see that partial fraction
decomposition helps in evaluating constant terms, and that only
part of the partial fraction is needed.

\section{Algorithm by Partial Fraction Decomposition}
Working in the field of iterated Laurent series has two
advantages. First, the expansion of a rational function into
Laurent series is determined by the total ordering ``$\preceq$ "
on its monomials, so we can temporarily forget its expansion as
long as we work in this field. Second, the fact that $F$ is a
rational function in $\lambda$ with  coefficients in a certain
\emph{field} permits us to apply the theory of partial fraction
decompositions.

Note that the idea of using partial fraction decompositions in
this context was first adopted by Stanley in \cite[p.
229--231]{stanley-rec}, but without the use of computers,  this
idea was thought to be impractical.

MacMahon's partition analysis always works in a ring like
$K[\Lambda, \Lambda^{-1}][[\mb{x}]]$, where $\Lambda^{-1}$ is
short for $\lambda^{-1}_1,\dots ,\lambda^{-1}_r$. This ring can be
embedded into a field of iterated Laurent series, such as $K\ll
\Lambda,\mb{x} \gg$.

While working in the field of iterated Laurent series, it is
convenient to use the operator $\pt_\lambda$, which is formally
defined by
$$\pt_\lambda \sum_{n=-\infty}^{\infty} a_n \lambda^n =\sum_{n=0}^\infty a_n \lambda^n,$$
whose validity is justified by the fundamental structure theorem.

MacMahon's operators can be realized as the following.
\begin{align}
\Oge F(\Lambda, \mb{x}) &=\left. \PT_{\Lambda} F(\Lambda, \mb{x})
\right|_{\Lambda= (1,\dots
,1)},\\
\Oeq F(\Lambda, \mb{x}) &=\ct_\Lambda F(\Lambda,\mb{x})=\left.
\PT_{\Lambda} F(\Lambda, \mb{x}) \right|_{\Lambda= (0,\dots ,0)}.
\end{align} So it suffices to find $\pt_\Lambda F$.
In fact, it is well-known that $\Oge$ can be realized by $\Oeq$ by
introducing new variables, just as the $\pt$ operators can be
realized by the $\ct$ operators (see \cite[Ch. 1]{xinthesis}). So
either an algorithm for $\pt_\Lambda F$ or an algorithm for
$\ct_\Lambda F$ will be sufficient for our purpose. Generally
speaking, $\pt$ is more suitable for the algorithm, and $\ct$ is
more suitable for theoretical analysis.

Now we need an algorithm to evaluate $\pt_\lambda F(\lambda)$ with
$$F(\lambda)= \frac{P(\lambda)}{\prod_{1\le i\le n} (\lambda^{j_i}-z_i)} $$
where $P(\lambda)$ is a polynomial in $\lambda$, $j_i$ are
nonnegative integers, and $z_i$ are independent of $\lambda$. Note
that we allow $z_i$ to be zero, so that the case of $P(\lambda)$
being a Laurent polynomial is covered. Our approach is different
from the previous algorithms, which deal with rational functions
expressed as in \eqref{e-2-mac-g} (the difference will be further
explained in the next section). It based on the following known
fact, which says that once the partial fraction decomposition of
$F$ is given, $\pt_\lambda F$ can be read off immediately.

\begin{thm}\label{t-ct-F}
Suppose that the factors in the denominator of $F$ are pairwise
relatively prime, and that the partial fraction decomposition of
$F$ is
$$F=f(\lambda)+\sum_{1\le i\le n} \frac{p_i(\lambda)}{\lambda^{j_i}-z_i},$$
where $f(\lambda)$ is a polynomial in $\lambda$, and
$p_i(\lambda)$ is a polynomial of degree less than $j_i$ for each
$i$. Then
\begin{align}\label{e-pt-F}
 \pt_\lambda F=f(\lambda) +\sum_{i}
\frac{p_i(\lambda)}{\lambda^{j_i}-z_i},
\end{align}
where the sum ranges over all $i$ such that $z_i \prec
\lambda^{j_i}$.
\end{thm}
\begin{proof}
The condition that $z_i$ is independent of $\lambda$ implies that
either $ \lambda^{j_i} \prec z_i$ or $z_i \prec \lambda^{j_i}$. In
the former case, we observe that the expansion of
$p_i(\lambda)/(\lambda^{j_i}-z_i)$ into Laurent series contains
only negative powers in $\lambda$, hence has no contribution when
applying $\pt_\lambda$. In the latter case, the expansion contains
only nonnegative powers in $\lambda$. Thus the the theorem
follows.
\end{proof}

To apply Theorem \ref{t-ct-F}, we need to know the partial
fraction decompositions of the given rational function. In fact,
we need only part of the partial fraction decompositions. Thus we
need an efficient algorithm for the partial fraction
decompositions. More ideally, an algorithm that only give us the
necessary parts. The classical algorithm does not seem to work
nicely. We use the new algorithm in \cite{xinparfrac} developed
from the following Theorem \ref{t-ppfraction-main1}.

To state the theorem, we need some concepts. Let $K$ be a field.
For $N,D\in K[t]$ with $D\ne 0$, $N/D$ can be uniquely written as
the summation of a polynomial $p$ and a proper fraction (or
rational function) $r/D$. We denote by $\poly(N/D)$ the polynomial
part, which is $p$, and by $\frr(N/D)$ the fractional part, which
is $r/D$.

Suppose that $N,D\in K[t]$ and $D$ is factored into pairwise
relatively prime factors $D=D_1\cdots D_k$. Then the ppfraction
(short for polynomial and proper fraction) expansion of $N/D$ with
respect to $D_1,\dots,D_k$ is the decomposition of $N/D$ as
$$N/D=p+r_1/D_1+\cdots r_k/D_k$$
such that $p,r_i$ are polynomials and $\deg(r_i)<\deg(D_i)$ for
every $i$. We denote the above $r_i/D_i$ by $\frr(N/D,D_i)$, the
fractional part of $N/D$ with respect to $D_i$.

\begin{thm}[Theorem 2.3 \cite{xinparfrac}]\label{t-ppfraction-main1}
For any $N,D\in K[t]$ with $D\ne 0$, if $D_1,\dots D_k \in K[t]$
are pairwise relatively prime, and $D=D_1\cdots D_k$, then
$$\frac{N}{D}= \poly\left(\frac{N}{D}\right) +\frr \left(\frac{N}{D}, D_1 \right)+\cdots +
\frr \left(\frac{N}{D}, D_k \right)$$ is the ppfraction expansion
of $N/D$ with respect to $(D_1,\dots ,D_k)$. Moreover, if
$1/(D_1D_i)= s_i/D_1+p_i/D_i$,  then
$$\frr(N/D,D_1)
=\frr(Ns_2s_3\cdots s_k/D_1).$$
\end{thm}

By Theorem \ref{t-ppfraction-main1}, we need two formulas to
develop our algorithm. One is  for the fractional part of
$p(\lambda)/(\lambda^j-a)$, and the other for  the partial
fraction decomposition of $(\lambda^j-a)^{-1}(\lambda^k-b)^{-1}$.
These are given as Propositions \ref{p-2-frac} and \ref{p-2-pfrac}
respectively.

Let $n \bmod {k}$ be the remainder of $n$ when divided by $k$. We
have

\begin{prop}\label{p-2-frac}
The fractional part of $p(\lambda)/(\lambda^j-a)$ can be obtained
by replacing $\lambda^d$ with $\lambda^{(d \bmod j)} a^{\lfloor
d/j \rfloor}$ in $p(\lambda)$ for all $d$, and dividing the result
by $\lambda^j-a$.
\end{prop}
\begin{proof}
By linearity, it suffice to show that the remainder of $\lambda^d$
when divided by $\lambda^j-a$
 equals
$\lambda^{(d \bmod j)} a^{\lfloor d/j \rfloor}$, which is trivial.
\end{proof}
It is easy to see that this operation takes time linear in the
number of nonzero terms of $p(\lambda)$, where we assume fast
arithmetic operations.
\begin{rem}\label{r-2-fracm}
Observe that the numerator of the fractional part of
$p(\lambda)/(\lambda^j-a)$ is always a Laurent polynomial in all
variables.
\end{rem}

\begin{lem}\label{l-2-pfrac}
For positive integers $j$ and $k$, if $a^k\ne b^j$, then the
following is a partial fraction expansion.
\begin{equation}
\frac{1}{(\lambda^j-a)(\lambda^k-b)} = \frac1{b^j-a^k}
\frr\left(\frac{\sum_{i=0}^{k-1} \lambda^{ij}a^{k-1-i}}
{\lambda^{k}-b}\right) -\frac1{b^j-a^k}
\frr\left(\frac{\sum_{i=0}^{j-1} \lambda^{ik}b^{j-1-i}}
{\lambda^{j}-a}\right)
\end{equation}
\end{lem}
\begin{proof}
First we show that if $a^k\ne b^j$, then $\lambda^j-a$ and
$\lambda^k-b$ are relatively prime. If not, say $\xi$ is their
common root in a field extension, then $\xi^j=a$ and $\xi^k=b$.
Thus we have $a^k=(\xi^j)^k=\xi^{jk}=(\xi^k)^j=b^j$, a
contradiction.

We have \begin{align*}
\frac{b^j-a^k}{(\lambda^j-a)(\lambda^k-b)}&=\frac{\lambda^{jk}-a^k}{(\lambda^j-a)(\lambda^k-b)}-
\frac{\lambda^{jk}-b^j}{(\lambda^j-a)(\lambda^k-b)}\\
&=\frac{\sum_{i=0}^{k-1} \lambda^{ij}a^{k-1-i}} {\lambda^{k}-b}
-\frac{\sum_{i=0}^{j-1} \lambda^{ik}b^{j-1-i}}{\lambda^j-a}.
\end{align*}
Now the polynomial part of
$\frac{b^j-a^k}{(\lambda^j-a)(\lambda^k-b)}$ is clearly $0$. Thus
the sum of the polynomial parts of the two terms on the right side
of the above equation also equals $0$. So taking the fractional
part of both sides and then dividing both sides by $b^j-a^k$ gives
the desired result.
\end{proof}

Now if $\gcd(j,k)$ is not $1$, then we can replace
$\lambda^{\gcd(j,k)}$ with $\mu$ and apply the above lemma. This
gives us the following result.

Let
$$\mathcal{F}(\lambda^j-a,\lambda^k-b)=\frac{\sum_{i=0}^{j'-1}
\lambda^{ik}b^{j'-1-i}}{a^{k'}-b^{j'}},$$ where $j'=j/\gcd(j,k)$
and $k'=k/\gcd(j,k)$.

\begin{prop}\label{p-2-pfrac}
For positive integers $j$ and $k$, if $a^k\ne b^j$, then we have
\begin{align}\label{e-2-pfracm}
\frr\left( \frac{1}{(\lambda^j-a)(\lambda^k-b)}, \lambda^j-a
\right)&= \frr \left(\frac{\mathcal{F}(\lambda^j-a,\lambda^k-b)
}{\lambda^{j}-a} \right),
\end{align}
\end{prop}
\begin{rem}
Note that a similar result appeared in \cite[Theorem 1]{george6},
but their proof was lengthy.
\end{rem}

Now by Theorem \ref{t-ppfraction-main1}, we have the following:
\begin{thm}\label{t-parfrac-F1}
With the notation of Theorem \ref{t-ct-F}, the polynomial
$p_s(\lambda)$ equals the remainder of
$$P(\lambda) \prod_{i=1,i\ne s}^n \mathcal{F}(\lambda^{j_s}-a_s, \lambda^{j_i}-a_i),$$ when divided by
$\lambda^{j_i}-z_i$ as a polynomial in $\lambda$.
\end{thm}

In Theorem \ref{t-ct-F}, we assumed that $\lambda^{j_i}-z_i$ and
$\lambda^{j_k}-z_k$ are relatively prime. Now let us consider the
case that $\lambda^{j_i}-z_i$ and $\lambda^{j_k}-z_k$ have a
nontrivial common factor. This happens if and only if
$z_i^{j_k}=z_k^{j_i}$, which can be easily checked. Andrews et al.
\cite{george6} suggested that we temporarily regard $z_i$ and
$z_j$ as two different variables. After the computation, we
replace them. We find an alternate approach, which has been
implemented in our computer program and will be discussed in the
next section.

Thus the above argument, Theorem \ref{t-ct-F}, and Theorem
\ref{t-parfrac-F1} together will give us an fast algorithm for
evaluating $\ct_\lambda F$.

\begin{rem}
From Remark \ref{r-2-fracm}, Theorem \ref{t-ct-F}, and Theorem
\ref{t-parfrac-F1}, we see that $\pt_{\lambda}F$ is
Elliott-rational when $F$ is. This is another way to prove Theorem
\ref{t-elliott}.
\end{rem}

\begin{exa}\label{ex-triangle}
Count all triples $(a,b,c)$ in $\NN^3$ such that they satisfy the
triangle inequalities.
\end{exa}
Similar problems have been done, such as counting non-congruent
triangles with integral side lengths \cite[Exercise 4.16]{EC1}. We
are going to illustrate our new approach by this example.
\begin{proof}[Solution] We solve the following three Diophantine inequalities:
$a+b-c\ge 0$, $b+c-a\ge 0$, and $c+a-b\ge 0$. The generating
function of these solutions is equal to $\Omega_\ge
F(\Lambda,\mb{x})$, where
$$F(\Lambda)=\frac{1}{((1-\lambda_1\lambda_3x_1/\lambda_2)(1-\lambda_1\lambda_2 x_2/\lambda_3)
(1-\lambda_2\lambda_3x_3/\lambda_1))}.$$

Although $F(\Lambda,\mb{x})$ is in
$K[\Lambda,\Lambda^{-1}][[\mb{x}]]$, we shall work in a field of
iterated Laurent series. We will chose $K\ll \Lambda,\mb{x}\gg$,
and $K\ll \Lambda, x_3,x_2,x_1\gg$ and compare the results.

We will apply $\Omega_\ge$ to $\lambda_3$, $\lambda_2$,
$\lambda_1$ subsequently. The first step, applying $\Omega_\ge$ to
$\lambda_3$ makes no difference for the two working fields.
Applying Theorems \ref{t-ct-F} and \ref{t-parfrac-F1} to the
factors of $F(\Lambda)$ containing $\lambda_3$, we get
\begin{multline*}
  \Omega_{\ge,\lambda_3}  F(\Lambda,\mb{x})
   =\\
   \frac{\lambda_2\lambda_1^{2}x_1}{ \left( {\lambda_1}^{2}x_1-\lambda_2^2x_3\right)  \left(
1-{\lambda_1}^{2}x_2x_1 \right)  \left( \lambda_1x_1-\lambda_2
\right) } -\frac {\lambda_1{\lambda_2}^{2}x_3}{
 \left( {\lambda_1}^{2}x_1-{\lambda_2}^{2}x_3
 \right)  \left( 1-{\lambda_2}^{2}x_2x_3 \right)  \left(
\lambda_1-\lambda_2x_3 \right)}
\end{multline*}
Denote by $F_1$ and $F_2$ the above two summands. At this stage,
we note that the expansion of $\left(
{\lambda_1}^{2}x_1-\lambda_2^2x_3\right)^{-1}$ does not exist in
$K[\Lambda,\Lambda^{-1}][[\mb{x}]]$, and generally there is no
advantage in getting rid of the factor
$\lambda_1^{2}x_1-\lambda_2^2x_3$ in the denominator by combining
the above two summands into one rational function. This will be
further explained in the next section.

Now when applying $\Omega_{\ge}$ on $\lambda_2$ to $F_1$ and
$F_2$, the results are different for the two working fields. Let
us look at $F_1$, especially the expansion of $\left(
{\lambda_1}^{2}x_1-\lambda_2^2x_3\right)^{-1}$. The expansion in
$K\ll \Lambda,\mb{x}\gg$ contains only nonnegative powers in
$\lambda_2$, while the expansion in $K\ll \Lambda, x_3,x_2,x_1\gg$
contains only negative powers in $\lambda_2$. The situation for
$F_2$ is similar. The conclusion is that working in $K\ll \Lambda,
x_3,x_2,x_1\gg$ is better: applying $\Omega_{\ge}$ on $\lambda_2$
to $F_1$ will gives us $0$, and we have

\begin{multline*}
    \Omega_{\ge,\lambda_3,\lambda_2} F(\Lambda,\mb{x})\\
    ={\frac {\lambda_{{1}} \left( \lambda_{{1}}+x_{{3}} \right) x_{{2}}}{
 \left( -x_{{3}}+{\lambda_{{1}}}^{2}x_{{2}} \right)  \left( -1+{
\lambda_{{1}}}^{2}x_{{2}}x_{{1}} \right)  \left( x_{{2}}x_{{3}}-1
 \right) }}+{\frac {\lambda_{{1}}x_{{3}}}{ \left( -x_{{3}}+{\lambda_{{
1}}}^{2}x_{{2}} \right)  \left( -1+x_{{1}}x_{{3}} \right)  \left(
-x_{ {3}}+\lambda_{{1}} \right) }}
\end{multline*}

Applying $\Omega_{\ge}$ on $\lambda_1$ to the two summands of the
above equation and simplifying gives us
\begin{align*}
\Omega_{\ge,\lambda_3,\lambda_2,\lambda_1} F(\Lambda,\mb{x})
={\frac {1+x_{{3}}x_{{2}}x_{{1}}}{ \left( 1-x_{{1}}x_{{3}} \right)
 \left( 1-x_{{2}}x_{{1}} \right)  \left(1- x_{{2}}x_{{3}} \right)
 }}.
\end{align*}
\end{proof}
\begin{rem}
  It is left to the reader to check that for the above example,
  $\CC\ll \Lambda ,x_2,x_3,x_1 \gg$ is the best working field.
\end{rem}

\section{The Maple Package}
\begin{lem}
Let $j_i$ be positive integers and let $z_i$ be monomials. If
$\lambda^{j_1}-z_1$ is not relatively prime to
$\lambda^{j_2}-z_2$, nor to $\lambda^{j_3}-z_3$, then
$\lambda^{j_2}-z_2$ and $\lambda^{j_3}-z_3$ are not relatively
prime.
\end{lem}
\begin{proof}
By the proof of Lemma \ref{l-2-pfrac}, we have
$z_1^{j_2}=z_2^{j_1}$ and $z_1^{j_3}=z_3^{j_1}$. It is easy to see
that $z_2^{j_3j_1}=z_3^{j_2j_1}$. Now the fact that $z_i$ is a
monomial (and hence has coefficient $1$) implies that
$z_2^{j_3}=z_3^{j_2}$.
\end{proof}

Thus to obtain the complete algorithm, we need to handle the
situation when $\lambda^{j_1}-z_1, \ldots,$ $ \lambda^{j_k}-z_k$
are not relatively prime to each other. For this situation, we
have not succeeded in applying the suggestion of the last section:
we tried to let $z_i=z_iv_i$ and do the computation, and finally
replace $v_i$ with $1$. But the problem is that the last step can
only be done after simplification, for which the rational function
will be too big for Maple to deal with. The following example
explains why this is not a fast approach: evaluating
$$\Omega_\ge \frac{1}{(1-\lambda x)^{10}(1-y/\lambda)^8}.$$

Our current program uses a modified ppfraction expansion as
follows. Suppose that $N$, $D$, $p_i$ belong to $K[\lambda]$, and
that $P=p_1\cdots p_k$ is relatively prime to $D$. Then we can
obtain a formula for $\frr(N/PD,P)$ satisfying our needs:

Write $N/(p_1D)=r_1/p_1+N_1/D$ with $\deg(r_1)<\deg(p_1)$. Then
$r_1/p_1=\frr(N/p_1D,p_1)$ can be easily obtained.

Now write $N_1/(p_2D)=r_2/p_2+N_2/D$. Then
$N/(p_1p_2D)=r_1/(p_1p_2)+r_2/p_2+N_2/D$.

In general, we have
$$\frac{N}{p_1\cdots p_kD}=\frac{r_1}{p_1\cdots p_k}+\cdots +\frac{r_k}{p_k}+\frac{N_k}{D},$$
with $\deg(r_i)<\deg(p_i)$. Now it is easy to see that
$$\frr\left(\frac{N}{PD},P\right)=\frac{r_1}{p_1\cdots p_k}+\cdots +\frac{r_k}{p_k}.$$
The recurrence formula for $r_i$ and $N_i$ is given by
$$ \frac{r_i}{p_i}=\frr\left(\frac{N_{i-1}}{p_iD},p_i\right) \quad \text{ and } N_i=\frac{N_{i-1}-r_iD}{p_i},$$
where $N_0=N$. Note that we shall let Maple compute $N_i$ with
respect to $\lambda$.

Now we can give the algorithm for computing $\pt_\lambda
F(\lambda)$ as follows.

\begin{enumerate}
    \item Collect the factors in the denominator of $F$ into
    several groups, such that the factors in different groups are
    relatively prime and factors in a same group are not.
    \item For each group having a contribution, find its corresponding fractional part of $F$.
    \item Take the sum of the results obtained from step 2, and add the
    polynomial part of $F$.
\end{enumerate}
\begin{rem}
We will simplify only if needed.

The factors in the denominator of $F$ that are independent of
$\lambda$ should be factored out to speed up the calculation. This
has been implemented in our computer program.
\end{rem}

The algorithm for $\Omega_\Lambda F(\Lambda,\mb{x})$ is described
as follows.

\begin{enumerate}
    \item Fix a total ordering on $\mb{x}$ and a total ordering on
    $\Lambda$. Suppose we are working in $\CC\ll \Lambda,\mb{x}
    \gg$.
    \item Eliminate $\lambda_1$ by computing $\pt_{\lambda_1} F$ and then replacing $\lambda_1$ with
    $1$.
    \item For each rational functions obtained from step 2,
    eliminate $\lambda_2$.
    \item Eliminate all the $\lambda$'s, and finally simplify.
\end{enumerate}

This approach partially solves the ``run-time explosion" problem
existing in Omega Calculus. Let us analyze a simple situation by
considering $\Omega_\ge F(\lambda)$, where
$$F(\lambda)=\frac{p(\lambda)}{\prod_{i=1}^m(1-x_i/\lambda)\prod_{j=1}^n (1-y_j\lambda)}.$$
The result after eliminating $\lambda$ and combining terms will
have a denominator of $mn$ factors: $(1-x_iy_j)$ with $1\le i\le
m$ and $1\le j\le n$. Such factors potentially contain the other
variables that are going to be eliminated. This explains the
existence of the run-time explosion problem.

In our approach, the result after eliminating $\lambda$ will be a
sum of $n$ rational functions (with a possible polynomial part),
each with a denominator of $m+n$ factors. Now it is crucial that
for each rational function, we can apply the theory of iterated
Laurent series to eliminate the other variables.

A Maple package implementing the above algorithm is available
online at \cite{xinEll}. Here is an example of how to use this
program after downloading this package. The current program uses
E\_{Oge}$(F,\mb{x},\Lambda)$ to compute $\Oge F(\Lambda,\mb{x})$
in the field $\CC\ll \Lambda, \mb{x}\gg$, where
 $\mb{x}$ is realized by $[x_1,\cdots ,x_n]$ in maple and
 $\Lambda$ is realized similarly.

\begin{exa}
Compute the generating function of $k$-gon partitions, which are
partitions that can be the side lengths of a $k$-gon.
\end{exa}
This problem was first studied in \cite{andrews9}, where the
generating functions of $k$-gon partitions are obtained only for
$k\le 6$ by using the authors' Omega package. We will discuss in
the next section about their formula for general $k$.

In the following $F(k)$ is the crude generating function of
$k$-gon partitions
$$F(k)= \frac{x_1a_1^{-1}}{(1-x_1\frac{a_k}{a_1})(1-x_2\frac{a_1a_k}{a_2})
\cdots
(1-x_{k-1}\frac{a_{k-2}a_k}{a_k-1})(1-x_k\frac{a_{k-1}}{a_k})},$$
where we use $a_i$ to replace $\lambda_i$. The function
$\text{test}(k)$ computes $\Oge F(k)$ and gives its normal
expression.

\begin{mapleinput}
\mapleinline{active}{1d}{read "Ell.mpl";}{}
\end{mapleinput}\begin{mapleinput}
\mapleinline{active}{1d}{F:=proc(k) }{}
\end{mapleinput}\begin{mapleinput}
\mapleinline{active}{1d}{product(1-q*a[k]*a[i-1]/a[i],i=2..k-1);}{}
\end{mapleinput}\begin{mapleinput}
\mapleinline{active}{1d}{q/a[1]/((1-q*a[k]/a[1])*
\end{mapleinput}\begin{mapleinput}
\mapleinline{active}{1d}{end:}{}
\end{mapleinput}
\begin{mapleinput}
\mapleinline{active}{1d}{va:=proc(k) seq(a[i],i=1..k) end:}{}
\end{mapleinput}\begin{mapleinput}
\mapleinline{active}{1d}{F(3);}{}
\end{mapleinput}
\mapleresult
\begin{maplelatex}
\mapleinline{inert}{2d}{q/(a[1]*(1-q*a[3]/a[1])*(1-q*a[3]*a[1]/a[2])*(1-q*a[2]/a[3]))}{%
\[\displaystyle
q{a_{{1}}}^{-1} \left( 1-{\frac {qa_{{3}}}{a_{{1}}}} \right) ^{-1}
\left( 1-{\frac {qa_{{3}}a_{{1}}}{a_{{2}}}} \right) ^{-1} \left(
1-{\frac {qa_{{2}}}{a_{{3}}}} \right) ^{-1}
\]
}

\end{maplelatex}
\begin{mapleinput}
\mapleinline{active}{1d}{E_Oge(
\end{mapleinput}
\mapleresult
\begin{maplelatex}
\mapleinline{inert}{2d}{-1/(q^3*(1/q^2-q^2)*(1/q^2-q)*(1-1/q^2))}{%
\[
-{q}^{-3} \left( {q}^{-2}-{q}^{2} \right) ^{-1} \left( {q}^{-2}-q
\right) ^{-1} \left( 1-{q}^{-2} \right) ^{-1}
\]
}

\end{maplelatex}
\begin{mapleinput}
\mapleinline{active}{1d}{test:=proc(n) F(n);va(n);E_Oge(
\end{mapleinput}\begin{mapleinput}
\mapleinline{active}{1d}{test(3);}{}
\end{mapleinput}
\mapleresult
\begin{maplelatex}
\mapleinline{inert}{2d}{-q^3/((q^4-1)*(q^3-1)*(q^2-1))}{%
\[\displaystyle
-{\frac {{q}^{3}}{ \left( {q}^{4}-1 \right)  \left( {q}^{3}-1
\right)  \left( {q}^{2}-1 \right) }}
\]
}

\end{maplelatex}
\begin{mapleinput}
\mapleinline{active}{1d}{test(4);}{}
\end{mapleinput}
\mapleresult
\begin{maplelatex}
\mapleinline{inert}{2d}{q^4*(q^3-q^2+1)/((q^3-1)*(-1+q)*(q^2-1)^2*(q^3+q^2+q+1)*(q^2-q+1))}{%
\[\displaystyle
{\frac {{q}^{4} \left( {q}^{3}-{q}^{2}+1 \right) }{ \left(
{q}^{3}-1 \right)  \left( -1+q \right)  \left( {q}^{2}-1 \right)
^{2} \left( {q}^{3}+{q}^{2}+q+1 \right)  \left( {q}^{2}-q+1
\right) }}
\]
}

\end{maplelatex}
\begin{mapleinput}
\mapleinline{active}{1d}{test(5);}{}
\end{mapleinput}
\mapleresult
\begin{maplelatex}
\mapleinline{inert}{2d}{-(q^10+q^9+q^8+q^7+q^6+q^5+q^4+q^3+q^2+q+1)*q^5/((q^3-1)*(-1+q)*(q^4-q^3+q-1)*(q^6+q^5-q-1)*(-1+q^8)*(q+1)*(q^2+1))}{%
\[\displaystyle
-{\frac { \left(
{q}^{10}+{q}^{9}+{q}^{8}+{q}^{7}+{q}^{6}+{q}^{5}+{q}^{4}+{q}^{3}+{q}^{2}+q+1
\right) {q}^{5}}{ \left( {q}^{3}-1 \right)  \left( -1+q \right)
\left( {q}^{4}-{q}^{3}+q-1 \right)  \left( {q}^{6}+{q}^{5}-q-1
\right)  \left( -1+{q}^{8} \right)  \left( q+1 \right)  \left(
{q}^{2}+1 \right) }}
\]
}

\end{maplelatex}
\begin{mapleinput}
\mapleinline{active}{1d}{test(6);}{}
\end{mapleinput}
\mapleresult
\begin{maplelatex}
\mapleinline{inert}{2d}{(q^12+q^11+q^10+q^9+q^8+2*q^7+q^6+q^5+q^3+q^2+q+1)*q^6/((q^3-1)*(q^5-1)*(q^8+q^6-q^2-1)*(q^4-1)*(q^2-1)*(q^6-q^4+q^2-1)*(q+1)*(q^4-q^3+q^2-q+1))}{%
\[
{\frac { \left(
{q}^{12}+{q}^{11}+{q}^{10}+{q}^{9}+{q}^{8}+2\,{q}^{7}+{q}^{6}+{q}^{5}+{q}^{3}+{q}^{2}+q+1
\right) {q}^{6}}{ \left( {q}^{3}-1 \right)  \left( {q}^{5}-1
\right)  \left( {q}^{8}+{q}^{6}-{q}^{2}-1 \right)  \left(
{q}^{4}-1 \right)  \left( {q}^{2}-1 \right)  \left(
{q}^{6}-{q}^{4}+{q}^{2}-1 \right)  \left( q+1 \right)  \left(
{q}^{4}-{q}^{3}+{q}^{2}-q+1 \right) }}
\]
}

\end{maplelatex}
\begin{mapleinput}
\mapleinline{active}{1d}{time(test(7));}{}
\end{mapleinput}
\mapleresult
\begin{maplelatex}
\mapleinline{inert}{2d}{34.765}{%
\[
 34.765
\]
}
\end{maplelatex}

All of the above are done in a personal computer. The time of
$\text{test}(7)$ is measured by seconds.

\begin{exa}
A Putnam problem (B3 on the 2000 Putnam examination) was
generalized in \cite{andrews7}. The generalized problems are
converted to evaluating $CP(T,k,c)=\Oge P(T,k,c)$ where
$$P(T,k,c)=\frac{1}{(1-x_1(a_1\cdots a_T)^ka_1^{-(k(T-1)+c)})
\cdots (1-x_T(a_1\cdots a_T)^k a_T^{-(k(T-1)+c)})}$$ for $k>c$ and
for $k<c$ we have
$$P(T,k,c)=\frac{1}{(1-x_1(a_1\cdots a_T)^{-k}a_1^{(k(T-1)+c)})
\cdots (1-x_T(a_1\cdots a_T)^{-k} a_T^{(k(T-1)+c)})}.$$ Note that
the case of $k=c$ is trivial.
\end{exa}

\begin{mapleinput}
\mapleinline{active}{1d}{read "Ell.mpl";}{}
\end{mapleinput}\begin{mapleinput}
\mapleinline{active}{1d}{P:=proc(T,k,c) if k>c then}{}
\end{mapleinput}\begin{mapleinput}
\mapleinline{active}{1d}{1/product(1-x[i]*product(a[j]^k,j=1..T)/a[i]^(k*(T-1)+c),i=1..T);}{}
\end{mapleinput}\begin{mapleinput}
\mapleinline{active}{1d}{else
1/product(1-x[i]/product(a[j]^k,j=1..T)*a[i]^(k*(T-1)+c),i=1..T)fi;
end:}{}
\end{mapleinput}\begin{mapleinput}
\mapleinline{active}{1d}{va:=proc(T) seq(a[i],i=1..T) end:}{}
\end{mapleinput}\begin{mapleinput}
\mapleinline{active}{1d}{vx:=proc(T) seq(x[i],i=1..T) end:}{}
\end{mapleinput}\begin{mapleinput}
\mapleinline{active}{1d}{CP:=proc(T,k,c) E_Oge(P(T,k,c),[vx(T)],[va(T)]);normal(
\end{mapleinput}

The case $T=3$, $k=2$ and $c=1$ is given explicitly.
\begin{mapleinput}
\mapleinline{active}{1d}{P(3,2,1);}{}
\end{mapleinput}
\mapleresult
\begin{maplelatex}
\mapleinline{inert}{2d}{1/((1-x[1]*a[2]^2*a[3]^2/a[1]^3)*(1-x[2]*a[1]^2*a[3]^2/a[2]^3)*(1-x[3]*a[1]^2*a[2]^2/a[3]^3))}{%
\[\displaystyle
 \left( 1-{\frac {x_{{1}}{a_{{2}}}^{2}{a_{{3}}}^{2}}{{a_{{1}}}^{3}}} \right) ^{-1}
 \left( 1-{\frac {x_{{2}}{a_{{1}}}^{2}{a_{{3}}}^{2}}{{a_{{2}}}^{3}}} \right) ^{-1}
 \left( 1-{\frac {x_{{3}}{a_{{1}}}^{2}{a_{{2}}}^{2}}{{a_{{3}}}^{3}}} \right) ^{-1}
\]
}

\end{maplelatex}
\begin{mapleinput}
\mapleinline{active}{1d}{CP(3,2,1);}{}
\end{mapleinput}
\mapleresult
\begin{maplelatex}
\mapleinline{inert}{2d}{-(x[1]^4*x[2]^4*x[3]^4+x[1]^3*x[2]^3*x[3]^3+x[2]^2*x[3]^2*x[1]^2+x[2]*x[1]*x[3]+1)/((-1+x[3]*x[2]^2*x[1]^2)*(-1+x[2]*x[3]^2*x[1]^2)*(x[2]^2*x[3]^2*x[1]-1))}{%
$\displaystyle -{\frac
{{x_{{1}}}^{4}{x_{{2}}}^{4}{x_{{3}}}^{4}+{x_{{1}}}^{3}{x_{{2}}}^{3}{x_{{3}}}^{3}+{x_{{2}}}^{2}{x_{{3}}}^{2}{x_{{1}}}^{2}+x_{{2}}x_{{1}}x_{{3}}+1}{
\left( -1+x_{{3}}{x_{{2}}}^{2}{x_{{1}}}^{2} \right)  \left(
-1+x_{{2}}{x_{{3}}}^{2}{x_{{1}}}^{2} \right)  \left(
{x_{{2}}}^{2}{x_{{3}}}^{2}x_{{1}}-1 \right) \mbox{}}}$}

\end{maplelatex}
\begin{mapleinput}
\mapleinline{active}{1d}{time(CP(3,3,1));}{}
\end{mapleinput}
\mapleresult
\begin{maplelatex}
\mapleinline{inert}{2d}{.171}{%
\[
 0.171
\]
}

\end{maplelatex}
\begin{mapleinput}
\mapleinline{active}{1d}{time(CP(3,2,3));}{}
\end{mapleinput}
\mapleresult
\begin{maplelatex}
\mapleinline{inert}{2d}{.391}{%
\[
 0.391
\]
}

\end{maplelatex}
\begin{mapleinput}
\mapleinline{active}{1d}{time(CP(3,1,3));}{}
\end{mapleinput}
\mapleresult
\begin{maplelatex}
\mapleinline{inert}{2d}{.235}{%
\[
 0.235
\]
}
\end{maplelatex}

\begin{mapleinput}
\mapleinline{active}{1d}{time(CP(3,1,4));}{}
\end{mapleinput}
\mapleresult
\begin{maplelatex}
\mapleinline{inert}{2d}{.686}{%
\[
 0.686
\]
}

\end{maplelatex}
\begin{mapleinput}
\mapleinline{active}{1d}{time(CP(3,1,5));}{}
\end{mapleinput}
\mapleresult
\begin{maplelatex}
\mapleinline{inert}{2d}{2.172}{%
\[
 2.172
\]
}

\end{maplelatex}
\begin{mapleinput}
\mapleinline{active}{1d}{CP(4,1,3):factor(
\end{mapleinput}
Maple will give us the following result:
\begin{multline*}
    {\frac { \left( {x_{{4}}}^{2}{x_{{3}}}^{2}{x_{{2}}}^{2}{x_{{1}}}^{2}+1
+x_{{2}}x_{{4}}x_{{3}}x_{{1}} \right)}{ \left(
x_{{4}}x_{{3}}{x_{{1}}}^{3}x_{{2}}-1 \right)
 \left( x_{{1}}x_{{4}}{x_{{3}}}^{3}x_{{2}}-1 \right)  \left( x_{{1}}{x
_{{4}}}^{3}x_{{3}}x_{{2}}-1 \right)  \left(
{x_{{2}}}^{3}x_{{4}}x_{{3} }x_{{1}}-1 \right) }} \\
 \big(
{x_{{2}}}^{3}{x_{{4}}}^{
3}{x_{{3}}}^{2}{x_{{1}}}^{3}+{x_{{2}}}^{3}{x_{{4}}}^{3}{x_{{3}}}^{3}{x
_{{1}}}^{3}+{x_{{4}}}^{3}{x_{{3}}}^{3}{x_{{2}}}^{2}{x_{{1}}}^{3}+
{x_{{
2}}}^{3}{x_{{4}}}^{2}{x_{{3}}}^{3}{x_{{1}}}^{3}+{x_{{2}}}^{3}{x_{{4}}}
^{3}{x_{{3}}}^{3}{x_{{1}}}^{2}\\+{x_{{1}}}^{2}x_{{4}}{x_{{3}}}^{2}x_{{2}
}+{x_{{1}}}^{2}{x_{{4}}}^{2}x_{{3}}x_{{2}}+{x_{{2}}}^{2}x_{{4}}x_{{3}}
{x_{{1}}}^{2}+x_{{4}}x_{{3}}x_{{2}}{x_{{1}}}^{2}+x_{{1}}{x_{{4}}}^{2}{
x_{{3}}}^{2}x_{{2}}\\+{x_{{2}}}^{2}x_{{4}}{x_{{3}}}^{2}x_{{1}}+x_{{4}}{x
_{{3}}}^{2}x_{{2}}x_{{1}}+{x_{{2}}}^{2}{x_{{4}}}^{2}x_{{3}}x_{{1}}+{x_
{{4}}}^{2}x_{{3}}x_{{2}}x_{{1}}+{x_{{2}}}^{2}x_{{4}}x_{{3}}x_{{1}}+1
 \big)
\end{multline*}

It will take Maple more than $5$ minutes to evaluate $CP(4,2,3)$.

We give a detailed comparison of the new algorithm and the Omega
package in Table \ref{table}, where the unit of the run-time is
seconds. Note that programs are not running on the same computer,
and that the data for the Omega package comes from
\cite{andrews7}.

\begin{table}
  \centering
  \caption{Comparison of our new method and the Omega package}\label{table}
\begin{tabular}{|c|c|c|c|c|c|c|}\hline
  run-time for CP& $(3,3,1)$ & $(3,2,3)$ & $(3,1,3)$ & $(3,1,4)$ & $(3,1,5)$ & $(4,1,3)$
  \\ \hline
 New method & 0.171 & 0.391 & 0.235 & 0.686 & 2.172 & 14.140 \\
 \hline
  Omega package & 7.14 & 58.95 & 100.55 & 643.86 & - & - \\ \hline
\end{tabular}
\end{table}

\begin{exa}
Magic squares of order $n$ are $n$ by $n$ matrices with integral
entries such that all the row sums and column sums are equal.
\end{exa}
The crude generating function for Magic square is:
\begin{align*}
\frac{1}{1-t(\lambda_1\cdots \lambda_n\mu_1\cdots
\mu_n)}\prod_{i=1}^n\prod_{j=1}^n
\frac{1}{1-x_{i,j}/(\lambda_i\mu_j)}.
\end{align*}
When evaluating the constant term in $\Lambda$ and $\mu$'s, we use
E\_Oeq instead of E\_Oge. Our maple program will reproduce the
result for $n=3$ quickly. For $n=4$, we get a sum of $96$ simple
rational functions, which is less than the $256$ of the Omega
package. Moreover, if we set $x_{i,j}=x$ at the beginning, and
finally replace $x$ with $1$, our program will reproduce the
formula for the case $n=5$ in about two minutes, which is the
generating function of the row sums for magic squares of order $5$
\cite[p.\ 234]{EC1}.

\section{Ways to Accelerate the Program}
Our program should have been accelerated by several ways, which
are not implemented due to the author's lack of programming
skills. These ways are list as follows and explained by examples.
we will manage to reduce of the number of rational functions of
the output, since the simplification of a sum of many rational
functions is a bottle neck for Maple (also Mathematica).

\begin{enumerate}

    \item The order of the variables to eliminate can make a difference
    for the computational time.

    \item The total ordering on the $x$'s can make a difference for the computational time,
    as will the total ordering on the $\lambda$'s.

\item The following alternative formula of \eqref{e-pt-F} can simplify the computation:
\begin{align}\label{e-pt-Fal}
 \pt_\lambda F(\lambda)=F(\lambda) -\sum_{i}
\frac{p_i(\lambda)}{\lambda^{j_i}-z_i},
\end{align}
where the sum ranges over all $i$ such that $z_i \succ
\lambda^{j_i}$.
\end{enumerate}

(1) is a well-known fact. To take advantage from it, we use the
fact that the number of rational functions produced by eliminating
$\lambda_i$ is equal to the number $cf(\lambda_i)$ of factors in
the denominator of $F$ that have contributions with respect to
$\lambda_i$. If $cf(\lambda_{i_0})$ is the smallest among all the
$cf(\lambda_i)$, then we shall eliminate the $\lambda_{i_0}$
first. Note that This way does not guarantee the best result.

The first part of (2) can be explained by Example
\ref{ex-triangle}, which gives a simple example of how to take
advantage of (2). The exact description will take time. The second
part is similar \cite[Example 2.5.13]{xinthesis}.

Using (3) might produce fewer rational functions. This happens
when the denominator of $F$ has more factors with contribution
than those factors without contribution. See the following
example.

\begin{exa}
Count all $k$-gons with nonnegative integral side lengths, which
are not required to be in an increasing order.
\end{exa}
\begin{proof}[Solution]
Suppose the side lengths of a $k$-gon is given by $a_1,\dots
,a_k$. Then we have $k$ inequalities, $a_1+\cdots +a_k\ge 2a_i$
for all $i$.

Using formula \eqref{e-pt-Fal} we can compute the generating
function of $k$-gons without computer. The eliminating order is
$\lambda_k,\dots ,\lambda_1$.
\begin{multline*}
\sum_{a_1+\cdots+a_k\ge 2a_i \text{ for all } i} x_1^{a_1}\cdots
x_k^{a_k}
 =
\Omega_\ge \frac{1}{(1-x_1\lambda_1\cdots
\lambda_k/\lambda_1^2)\cdots (1-x_k \lambda_1\cdots \lambda_k
/\lambda_k^2)}\\
=\Omega_\ge \frac{1}{(1-x_1\lambda_1\cdots
\lambda_{k-1}/\lambda_1^2)\cdots (1-x_{k-1} \lambda_1\cdots
\lambda_{k-1}
/\lambda_{k-1}^2)(1-x_k\lambda_1\cdots\lambda_{k-1})}\\
-\frac{x_k\lambda_1\cdots
\lambda_{k-1}}{(1-x_1x_k\lambda_1^2\cdots
\lambda_{k-1}^2/\lambda_1^2)\cdots (1-x_{k-1}x_k \lambda^2_1\cdots
\lambda_{k-1}^2 /\lambda_{k-1}^2)(1-x_k\lambda_1\cdots
\lambda_{k-1})}.
\end{multline*}
Now notice that $\Omega_\ge $ acting on the second term is simply
obtained by replacing $\lambda_i$ with $1$. Repeat the above
computation, we get the final generating function:
$$\frac{1}{(1-x_1)\cdots (1-x_k)}-\sum_{i=1}^k \frac{x_i}{(1-x_1x_i)\cdots (1-x_kx_i)}. $$
\end{proof}

However, it is probably better not to use \eqref{e-pt-Fal} if the
total degree of those factors without a contribution is much
greater than those factors with a contribution. Also note that
this formula is not easy to apply for the CT operator, we shall
use $\ct_\lambda F(\lambda)=\ct_\lambda F(1/\lambda)$ instead.

\vspace{3mm} There are also ways that may speed up the
computation, but are not easy to implement by the computer. The
following example is simplified by using different parameters for
a given problem.

\begin{exa}
Count $k$-gon partitions (revisited).
\end{exa}
An exact formula for the generating function of $k$-gon partitions
was given in \cite[Theorem 1]{andrews9}. Here we give a simple
proof by using different parameters and formula \eqref{e-pt-Fal}.
\begin{proof}[Solution]
The problem is to find all $(a_1,\dots,a_k)\in \PP^k$ such that
$1\le a_1\le a_2\le \cdots \le a_k$, and $a_1+\cdots +a_{k-1}>
a_k$.

Let $b_1=a_1-1$, $b_2=a_2-a_1$,\dots, $b_k=a_k-a_{k-1}$. Then
$a_i=1+b_1+\cdots +b_i$ for all $i$, and it suffices to find all
$b_i$ such that $b_i\ge 0$, and $k-3+(k-2)b_1+(k-3)b_2+\cdots
b_{k-2}\ge b_k$. Thus the generating function for these $b_i$ are
given by
\begin{multline*}
\Oge \frac{\lambda^{k-3}}{(1-x_1\lambda^{k-2})\cdots
(1-x_{k-2}\lambda)(1-x_{k-1})(1-x_k/\lambda)}\\
=\frac{1}{(1-x_1)\cdots
(1-x_k)}-\frac{x_k^{k-2}}{(1-x_1x_k^{k-2})(1-x_2x_k^{k-3})\cdots
(1-x_{k-1})(1-x_k)}.
\end{multline*}
Now it is easy to convert this formula to \cite[Theorem
1]{andrews9}.
\end{proof}

\vspace{3mm} {\bf Acknowledgment.} The author is very grateful to
his advisor Ira Gessel.

\bibliographystyle{amsplain}
\providecommand{\bysame}{\leavevmode\hbox
to3em{\hrulefill}\thinspace}
\providecommand{\MR}{\relax\ifhmode\unskip\space\fi MR }
\providecommand{\MRhref}[2]{%
  \href{http://www.ams.org/mathscinet-getitem?mr=#1}{#2}
} \providecommand{\href}[2]{#2}

\end{document}